\def\ver{Apr. 14, 2004, v.5}
\documentstyle{amsppt}
\magnification=1200
\hsize=6.5truein
\vsize=8.9truein
\topmatter
\title Multiplier Ideals,
$ V $-Filtration, and Spectrum
\endtitle
\author Nero Budur and Morihiko Saito
\endauthor
\abstract
For an effective divisor on a smooth algebraic variety or a
complex manifold, we show that the associated multiplier ideals
coincide essentially with the filtration induced by the filtration V
constructed by B.~Malgrange and M.~Kashiwara.
This implies another proof of a theorem of L.~Ein, R.~Lazarsfeld,
K.E.~Smith and D.~Varolin that any jumping coefficient in the
interval (0,1] is a root of the Bernstein-Sato polynomial
up to sign.
We also give a refinement (using mixed Hodge modules) of the
formula for the coefficients of the spectrum
for exponents not greater than one or greater than the dimension
of the variety minus one.
\endabstract
\endtopmatter
\tolerance=1000
\baselineskip=12pt
\def\ssbull{\raise.2ex\hbox{${\scriptscriptstyle\bullet}$}}
\def\msum{\hbox{$\sum$}}
\def\mopls{\hbox{$\bigoplus$}}
\def\mcup{\hbox{$\bigcup$}}
\def\mcap{\hbox{$\bigcap$}}
\def\bC{{\Bbb C}}
\def\be{\bold{e}}
\def\bD{{\Bbb D}}
\def\bN{{\Bbb N}}
\def\bQ{{\Bbb Q}}
\def\bZ{{\Bbb Z}}
\def\boR{\bold{R}}
\def\cB{{\Cal B}}
\def\cD{{\Cal D}}
\def\cG{{\Cal G}}
\def\cJ{{\Cal J}}
\def\cO{{\Cal O}}
\def\tD{{\widetilde{D}}}
\def\tH{{\widetilde{H}}}
\def\tom{{\widetilde{\omega}}}
\def\Sp{\hbox{\rm Sp}}
\def\Ker{\hbox{\rm Ker}}
\def\supp{\hbox{\rm supp}\,}
\def\codim{\hbox{\rm codim}\,}
\def\Gr{\hbox{\rm Gr}}
\def\can{\text{\rm can}}
\def\red{\text{\rm red}}
\def\simto{\buildrel\sim\over\longrightarrow}
\def\SameAuthor{\vrule height3pt depth-2.5pt width1cm}

\document\noindent
\centerline{\bf Introduction}\footnote""{{\it Date}\,: \ver}

\bigskip\noindent
Let
$ X $ be a smooth complex algebraic variety or a complex
manifold, and
$ D $ be an effective divisor on
$ X $ with a defining equation
$ f $.
The multiplier ideal
$ \cJ(D) $ is a coherent sheaf of ideals of the structure
sheaf
$ \cO_{X} $,
and can be defined by using an embedded resolution
$ \pi : (X',D') \to (X,D) $,
see [8], [10], [17], [21].
This is defined also for the
$ \bQ $-divisors
$ \alpha D $ with
$ \alpha > 0 $, and we get a decreasing family
$ \{\cJ(\alpha D)\}_{\alpha\in\bQ} $, where
$ \cJ(\alpha D) = \cO_{X} $ for
$ \alpha \le 0 $.
By construction there exist positive rational numbers
$ 0 < \alpha_{1} < \alpha_{2} < \cdots $ such that
$ \cJ(\alpha_{j}D) = \cJ(\alpha D) \ne \cJ(\alpha_{j+1}D) $
for
$ \alpha_{j} \le \alpha < \alpha_{j+1} $ where
$ \alpha_{0} = 0 $,
see (1.1).
These numbers
$ \alpha_{j} \,(j > 0) $ are called the jumping coefficients
(or numbers) of the multiplier ideals associated to
$ D $.
We define the graded pieces
$ \cG(D,\alpha) $ to be
$ \cJ((\alpha -\varepsilon)D)/\cJ(\alpha D) $ for
$ 0 < \varepsilon \ll 1 $.

Let
$ i_{f} : X \to Y := X \times \bC $ denote the embedding
by the graph of
$ f $,
and
$ t $ be the coordinate of
$ \bC $.
Let
$ \cB_{f} \,(= \cO_{X}\otimes_{\bC}\bC[\partial_{t}]) $
denote the direct image of the left
$ \cD_{X} $-module
$ \cO_{X} $ by
$ i_{f} $.
B.~Malgrange [20] constructed the filtration
$ V $ on
$ \cB_{f} $ and M.~Kashiwara [15] did it in a more
general case, see also [16].
We index
$ V $ decreasingly by rational numbers so that the action of
$ \partial_{t}t - \alpha $ on
$ \Gr_{V}^{\alpha}\cB_{f}\,(= V^{\alpha}\cB_{f}/
V^{>\alpha}\cB_{f}) $ is locally nilpotent.
We denote also by
$ V $ the induced filtration on
$ \cO_{X} \,(= \cO_{X}\otimes 1) $.
Then

\medskip\noindent
{\bf 0.1.~Theorem.}
{\it
We have
$ V^{\alpha}\cO_{X} = \cJ(\alpha D) $ if
$ \alpha $ is not a jumping coefficient.
In general,
$ \cJ(\alpha D) = V^{\alpha + \varepsilon}\cO_{X} $ and
$ V^{\alpha}\cO_{X} = \cJ((\alpha -\varepsilon)D) $ for any
$ \alpha \in \bQ $ if
$ \varepsilon > 0 $ is sufficiently small.
}

\medskip
This is actually an immediate consequence of [24], 3.5
and [23], 3.3.17.
Indeed, the normal crossing case was proved in the former (see
also (2.3) below) and the general case follows from it using
the theory of bifiltered direct images developed in the latter,
see (3.2) and (3.4) below.
By Theorem (0.1) we have

\medskip\noindent
{\bf 0.2.~Corollary.}
{\it
$ \Gr_{V}^{\alpha}\cO_{X} = \cG(D,\alpha) $ for any
$ \alpha \in \bQ $.
}

\medskip
This corollary together with Theorem (0.3) below has
been conjectured in [3], and
proved in [5] using [24], 3.5 together with [24], 2.14
instead of [23], 3.3.17.
As a corollary of (0.2), we get another proof of a theorem of
L.~Ein, R.~Lazarsfeld, K.E.~Smith and D.~Varolin [12] that any
jumping coefficient in the interval
$ (0,1] $ is a root of the
$ b $-function (i.e. the Bernstein-Sato polynomial) up to
sign.

If
$ D $ is reduced, let
$ \rho : \tD \to D $ be a resolution of singularities.
It is well known that
$ \rho_{*}\omega_{\tD} \subset \omega_{D} $ is
independent of
$ \tD $ (see [13]).
Let
$ \tom_{D} = \rho_{*}\omega_{\tD} $.
Combining (0.2) with [25], we get

\medskip\noindent
{\bf 0.3.~Theorem.}
{\it If
$ D $ is reduced, there exists a decreasing filtration
$ V $ on
$ \omega_{D}/\tom_{D} $ indexed by
$ (0,1] \cap \bQ $ {\rm (}i.e.
$ V^{>0}(\omega_{D}/\tom_{D}) = \omega_{D}/\tom_{D},
V^{>1}(\omega_{D}/\tom_{D}) = 0) $ such that
$$
\Gr_{V}^{\alpha}(\omega_{D}/\tom_{D}) =
\omega_{D}\otimes \cG(D,\alpha)\quad
\text{for}\,\, 0 < \alpha < 1,
$$
and
$ \Gr_{V}^{1}(\omega_{D}/\tom_{D}) $ is a
quotient of
$ \omega_{D}\otimes \cG(D,1) $.
}

\medskip
This means that the
$ \cG(D,\alpha) $ measure the non rationality of the
singularity of
$ D $.
In the isolated singularity case, (0.3) is related to [18], [19],
[29].
For
$ \alpha = 1 $ we have a more precise formula, see (3.5--6) below.

Let
$ x \in X $.
J.~Steenbrink ([26], [27]) defined the spectrum
$ \Sp(f,x) $ by using the monodromy and the mixed Hodge
structure on the vanishing cohomology groups
$ H^{j}(F_{x},\bC) $, where
$ F_{x} $ denotes the Milnor fiber of
$ f $ around
$ x $.
In this paper we use the normalization such that the
spectrum is a fractional polynomial
$ \msum_{0<\alpha <n} n_{\alpha}t^{\alpha} $ where
$ n = \dim X $,
see (5.2).
In [3] and [4], the ``localized'' graded pieces
$ \cG(D,x,\alpha) $, which are quotients of
$ \cG(D,\alpha) $, were defined.
Using a formula for nearby cycles in [9] together with the
theory of cyclic covers and Serre duality,
the following theorem was proved:

\medskip\noindent
{\bf 0.4.~Theorem} [4].
{\it
$ n_{\alpha} = \dim \cG(D,x,\alpha) $ for
$ 0 < \alpha \le 1 $.
}

\medskip
By construction there is a coherent sheaf
$ K_{\alpha} $ on the embedded resolution
$ (X',D') $ of
$ (X,D) $ such that
$ \cG(D,x,\alpha) $ is the global sections of
$ K_{\alpha} $.
Here
$ E := \pi^{-1}(x)_{\red} $ is assumed to be a divisor
with normal crossings, and
$ K_{\alpha} $ is actually the restriction of
$ \cG(D,x,\alpha) $ to
$ E $ as an
$ \cO $-module
(up to the twist by the relative dualizing sheaf).
We have also
$ K'_{\alpha} $ by twisting
$ K_{\alpha} $ with a certain line bundle, see (5.4.1).
We prove that these sheaves coincide with the first nonzero
piece of the Hodge filtration of the mixed Hodge modules which
calculate the vanishing cohomology groups with or without
compact supports, see (4.3).
This implies

\medskip\noindent
{\bf 0.5.~Theorem.}
{\it For
$ j \in \bZ $ and
$ \alpha \in (0,1] $, there are natural isomorphisms
$$
\aligned
H^{j}(E, K_{\alpha})
&= F^{n-1}H^{j+n-1}(F_{x},\bC)_{\be(-\alpha)},
\\
H^{j}(E, K'_{\alpha})
&= F^{n-1}H_{c}^{j+n-1}(F_{x},\bC)_{\be(-\alpha)}.
\endaligned
$$
where
$ H^{j}(F_{x},\bC)_{\lambda} $ is the
$ \lambda $-eigenspace of the monodromy and
$ \be(-\alpha) = \exp(-2\pi i\alpha) $.
}

\medskip
This gives a more direct proof of (0.4).
Using duality, (0.5) implies

\medskip\noindent
{\bf 0.6.~Corollary.}
{\it
$ n_{\alpha} = \chi (E,K'_{n-\alpha}) $ for
$ n - 1 < \alpha < n $.
}

\medskip
This was also proved in [3].
Note that (0.5) is useful for an explicit calculation of
the spectrum.
In the case
$ n = 3 $, it is reduced to that of
$ \dim \cG(D,x,\alpha), \chi (E, K'_{3-\alpha}) $ and the
monodromy.

The first author would like to thank Professor L.~Ein for
several discussions and comments on this work.
The second author would like to thank Professor S.~Mukai
for drawing his attention to the first author's work [4].

In Sect.~1 we recall the definition of multiplier ideal.
In Sect.~2 we review the theory of
$ V $-filtrations, and calculate it in the normal crossing
case.
In Sect.~3 we review the theory of bifiltered direct images
and prove Theorem (0.1).
In Sect.~4 we study the localizations of nearby cycle sheaves
for the proof of (0.5).
In Sect.~5 we recall the definition of spectrum and prove
Theorem (0.5).

\bigskip\bigskip\centerline{{\bf 1. Multiplier Ideals}}

\bigskip\noindent
{\bf 1.1.}
Let
$ X $ be a smooth complex algebraic variety or a complex
manifold, and
$ D $ be an effective divisor on it.
Let
$ \pi : (X',D') \to (X,D) $ be an embedded resolution, i.e.
$ D' := \pi^{*}D $ is a divisor with normal crossings.
Here
$ \pi $ is assumed to be projective.
Let
$ D'_{i} $ be the irreducible components of
$ D' $ with
$ m_{i} $ the multiplicity so that
$ D' = \msum_{i} m_{i}D'_{i} $.
For positive rational numbers
$ \alpha $,
we define
$$
\cJ(\alpha D) = \pi_{\ssbull}(\omega_{X'/X}\otimes
\cO_{X'}(-\msum_{i} [\alpha m_{i}]D'_{i})),
\leqno(1.1.1)
$$
and
$ \cJ(\alpha D) = \cO_{X} $ for
$ \alpha \le 0 $.
Here
$ \omega_{X'/X} $ is the relative dualizing sheaf
$ \omega_{X'}\otimes \pi^{*}\omega_{X}^{\vee} $,
and
$ \pi_{\ssbull} $ denotes the sheaf theoretic direct image.
They form a decreasing family of coherent subsheaves of
$ \cO_{X} $,
and there exist (at least locally) a strictly increasing
sequence of positive rational numbers
$ \{\alpha_{j}\}_{j>0} $ (called the jumping coefficients)
such that
$$
\cJ(\alpha_{j}D) = \cJ(\alpha D) \ne \cJ(\alpha_{j+1}D)
\quad \text{for}\,\,
\alpha_{j} \le \alpha < \alpha_{j+1},
\leqno(1.1.2)
$$
where
$ \alpha_{0} = 0 $.
We define the graded pieces by
$ \cG(D,\alpha) = \cJ((\alpha -\varepsilon)D)/
\cJ(\alpha D) $ for
$ 0 < \varepsilon \ll 1 $.
It is well-known (see [8], [10], [21]) that
the higher direct images vanish, i.e.
$$
R^{j}\pi_{\ssbull}(\omega_{X'/X}\otimes
\cO_{X'}(-\msum_{i} [\alpha m_{i}]D'_{i})) = 0
\quad\text{for}\,\, j > 0.
\leqno(1.1.3)
$$
So the
$ \cG(D,\alpha) $ are isomorphic to the direct image of
the quotient sheaves.

\medskip\noindent
{\bf 1.2.~Remark.}
Originally, the multiplier ideal
$ \cJ(\alpha D) $ is defined by the local integrability of
$ |g|^{2}/|f|^{2\alpha} $ for
$ g \in \cO_{X} $, see [8], [21].
This interpretation justifies the definition of the graded
pieces
$ \cG(D,\alpha) $.

\newpage
\centerline{{\bf 2. $ V $-Filtration}}

\bigskip\noindent
{\bf 2.1.}
Let
$ X $ be a smooth complex algebraic variety or
a complex manifold of dimension
$ n $,
and
$ D $ be an effective divisor with a defining equation
$ f $.
Let
$$
\cB_{f} = \cO_{X}\otimes_{\bC}\bC[\partial_{t}]
\leqno(2.1.1)
$$
as in the introduction.
The action of
$ \cD_{Y} $ on
$ \cO_{X}\otimes_{\bC}\bC[\partial_{t}] $ is given by
$$
\alignat 2
g(h\otimes {\partial}_{t}^{j})
&= gh\otimes {\partial}_{t}^{j},
&\quad \xi (h\otimes {\partial}_{t}^{j})
&= \xi h\otimes {\partial}_{t}^{j} - (\xi f)h\otimes
{\partial}_{t}^{j+1}
\\
\partial_{t}(h\otimes {\partial}_{t}^{j})
&= h\otimes{\partial}_{t}^{j+1},
&\quad t(h\otimes {\partial}_{t}^{j})
&= fh\otimes {\partial}_{t}^{j} - jh\otimes
{\partial}_{t}^{j-1}
\endalignat
$$
for
$ g, h \in \cO_{X}, \zeta \in \Theta_{X} $, i.e.
$ \cO_{X}\otimes_{\bC}\bC[\partial_{t}] $ is identified with
$ \cO_{X}\otimes_{\bC}\bC[\partial_{t}]\delta (f-t) $ where
$ \delta(f-t) $ is the delta function.

The filtration
$ V $ on
$ \cB_{f} $ (see [15], [20]) is an exhaustive decreasing
filtration of coherent
$ \cD_{X} $-submodules, and is characterized by the following
properties:

\medskip\noindent
(i)
$ t(V^{\alpha}\cB_{f}) \subset V^{\alpha +1}\cB_{f} $,
$ \partial_{t}(V^{\alpha}\cB_{f}) \subset V^{\alpha -1}
\cB_{f} $ for any
$ \alpha \in \bQ $,

\medskip\noindent
(ii)
$ t(V^{\alpha}\cB_{f}) = V^{\alpha +1}\cB_{f} $ for
$ \alpha > 0 $,

\medskip\noindent
(iii) the action of
$ \partial_{t}t - \alpha $ on
$ \Gr_{V}^{\alpha}\cB_{f} $ is locally nilpotent.

\medskip\noindent
Here
$ \Gr_{V}^{\alpha} = V^{\alpha}/V^{>\alpha} $ with
$ V^{>\alpha} = \mcup_{\beta>\alpha}V^{\beta} $, and
$ V $ is indexed discretely by
$ \bQ $ so that for some positive integer
$ m $
$$
V^{\alpha} = V^{j/m}\quad \text{for}\,\,
(j-1)/m < \alpha \le j/m.
\leqno(2.1.2)
$$
In this case
$ m $ can be taken to be the order of the semisimple part
of the monodromy on the nearby cycles.

The Hodge filtration
$ F $ on
$ \cB_{f} $ is defined by the order of
$ \partial_{t} $ shifted by
$ - n $,
i.e.
$$
F_{p-n}\cB_{f} = \mopls_{0\le j\le p} \cO_{X}\otimes
{\partial}_{t}^{j}.
\leqno(2.1.3)
$$
This shift comes from the shift of the Hodge filtration
$ F $ on
$ \cO_{X} $ which is defined by
$ \Gr_{p}^{F} = 0 $ for
$ p \ne - n $.
(Actually this is for the corresponding right
$ \cD $-module
$ \omega_{X} $.
We use this because the direct image by a closed embedding
requires a shift of filtration if we use the usual filtration
for left
$ \cD $-modules.)
We have
$$
F_{-n}\Gr_{V}^{\alpha}\cB_{f} = 0 \quad
\text{for}\,\, \alpha \le 0,
\leqno(2.1.4)
$$
because
$ \partial_{t} : \Gr_{V}^{\alpha+1}(\cB_{f},F) \to
\Gr_{V}^{\alpha}(\cB_{f},F[-1]) $ is strictly surjective for
$ \alpha \le 0 $, see [23], 3.2.1 and 5.1.4.
In particular
$ V^{>0}\cO_{X} = \cO_{X} $.
This is related to [14].

\medskip\noindent
{\bf 2.2.~Normal crossing case.}
Assume
$ D = f^{-1}(0) $ is a divisor with normal crossings.
Let
$ (x_{1}, \dots, x_{n}) $ be a local coordinate at
$ x \in X $ such that
$ f = \prod_{i} x_{i}^{m_{i}} $ (\'etale locally
in the algebraic case), where
$ m_{i} $ are nonnegative integers.
For a multi index
$ \nu = (\nu_{1}, \dots, \nu_{n}) \in \bN^{n} $,
let
$ x^{\nu} = \prod_{i} x_{i}^{{\nu}_{i}} $.
We define the filtration
$ V^{\alpha}\cO_{X,x} $ to be the
$ \cO_{X,x} $-submodule generated by
$$
x^{\nu}\quad \text{with}\,\, \nu_{i} + 1 \ge m_{i}\alpha \quad
\text{for any}\,\, i.
\leqno(2.2.1)
$$
This is generated by one element, and
$$
f(V^{\alpha}\cO_{X,x}) = V^{\alpha +1}\cO_{X,x} \quad
\text{for}\,\, \alpha > 0.
\leqno(2.2.2)
$$
By [24], 3.4,
$ V^{\alpha}\cB_{f,x} $ is generated over
$ \cD_{X,x} $ by
$$
V^{\alpha +j}\cO_{X,x}\otimes {\partial}_{t}^{j}\quad
\text{with}\,\, j \le \max\{1 - \alpha,0\}.
\leqno(2.2.3)
$$
Indeed, this follows from
$$
(\partial/\partial x_{i})x_{i}(x^{\nu}\otimes
{\partial}_{t}^{j}) = (\nu_{i} + 1 - m_{i}(s + j))
x^{\nu}\otimes {\partial}_{t}^{j},
\leqno(2.2.4)
$$
where
$ s = \partial_{t}t $.

For
$ \alpha \in (0,1] $,
let
$ D(\alpha) $ be the union of the irreducible components
$ D_{i} $ whose multiplicity
$ m_{i} $ satisfies
$ m_{i}\alpha \in \bZ $.
Then we have
$$
\supp \Gr_{V}^{\alpha}\cB_{f} \subset D(\alpha)
\leqno(2.2.5)
$$

\medskip\noindent
{\bf 2.3.~Proposition.}
{\it With the assumption of {\rm (2.2),} we have for any
$ \alpha $
}
$$
F_{-n}V^{\alpha}\cB_{f,x} = V^{\alpha}\cO_{X,x}
\otimes 1.
\leqno(2.3.1)
$$

\medskip\noindent
{\it Proof.}
See [24], Prop. 3.5.

\medskip\noindent
{\bf 2.4.~Remark.}
The equality (2.3.1) is not a simple corollary of (2.2.3),
and we need some calculation.
Note that (2.3.1) is equivalent to Theorem (0.1) in the
normal crossing case using
$$
\min\{r \in \bZ : r + 1 \ge \beta \} =
\max\{r \in \bZ : r < \beta \},
\leqno(2.4.1)
$$
where
$ \beta = m_{i}\alpha $.
To deduce the general case from this, we need the theory of
bifiltered direct images as explained in the next section.

\bigskip\bigskip\centerline{{\bf 3. Bifiltered Direct Images}}

\bigskip\noindent
{\bf 3.1.} Let
$ \pi : X' \to X $ be a projective morphism of complex
manifolds or a proper morphism of smooth complex algebraic
varieties.
Let
$ n = \dim X $,
$ m = \dim X' $,
and
$ f' = f\pi $.
Then we can define
$ \cB_{f'} $ on
$ Y' = X' \times \bC $ similarly.
We factorize
$ \pi $ by the composition of the embedding by the graph
$$
i_{\pi} : X' \to X'' := X' \times X
$$
and the projection
$ pr : X'' \to X $.
Let
$ (x_{1}, \dots, x_{n}) $ be a local coordinate system of
$ X $,
and put
$ \partial_{i} = \partial /\partial x_{i} $.
Let
$ h_{i} $ be the pull-back of
$ x_{i} $ by
$ pr $.
Then
$ \cB'' := (i_{\pi}\times id)_{*}\cB_{f'} $ is defined
locally on
$ X $ by
$$
\cB_{f'}\otimes_{\bC}\bC[\partial_{1}, \dots,
\partial_{n}] = \cO_{X'}\otimes_{\bC}\bC[\partial_{1},
\dots, \partial_{n}, \partial_{t}],
\leqno(3.1.1)
$$
and the action of a vector field
$ \xi $ on
$ X' $ is given by
$$
\xi(g \otimes \partial^{\nu}{\partial}_{t}^{j}) =
(\xi g)\otimes \partial^{\nu}{\partial}_{t}^{j} - \msum_{i}
(\xi h_{i})g \otimes \partial_{i}\partial^{\nu}
{\partial}_{t}^{j} -
(\xi f')g \otimes \partial^{\nu}{\partial}_{t}^{j+1}
\leqno(3.1.2)
$$
for
$ g \in \cO_{X'} $,
where
$ \nu = (\nu_{1}, \dots, \nu_{n}) $ is a multi index.
The filtration
$ F $ on
$ \cB'' $ is given by the total order of
$ \partial_{i}, \partial_{t} $,
and
$ V $ is generated over
$ \bC[\partial_{1}, \dots, \partial_{n}] $ by
$ V $ on
$ \cB_{f'} = \cO_{X'}\otimes_{\bC}\bC[\partial_{t}] $.
Globally we need the twist by the relative dualizing sheaf
[2].
For example
$$
F_{-m}\cB'' = \pi^{*}\omega_{X}^{\vee},
\leqno(3.1.3)
$$
where
$ \omega_{X}^{\vee} $ is the dual of the line bundle
$ \omega_{X} $.

The bifiltered direct image
$ (\pi\times id)_{*}(\cB_{f'};F,V) $ is then defined to be the
direct image of
$ \cB'' $ by
$ pr $.
It is defined by the sheaf theoretic direct image (using the
canonical flasque resolution) of the relative de Rham complex
$ DR_{X'\times X/X}(\cB'';F,V) $ whose
$ i $-th component is
$$
(\cB'';F[-i],V) \otimes_{\cO_{X'}} \Omega_{X'}^{i+m}.
$$
Here the shift of filtration for an increasing filtration is
given by
$ (F[m])_{p} = F_{p-m} $,
see [6].
For
$ p = -m $ we have by definition
$$
F_{-m}V^{\alpha}(\pi\times id)_{*}\cB_{f'} =
\boR\pi_{\ssbull}(\omega_{X'/X}\otimes V^{\alpha}
\cO_{X'}).
\leqno(3.1.4)
$$

\medskip\noindent
{\bf 3.2.~Proposition.}
{\it The direct image
$ (\pi\times id)_{*}(\cB_{f'};F,V) $ in {\rm (3.1)} is bistrict
{\rm (}i.e. strict in the sense of {\rm [23]),} and the induced
filtration
$ V $ on the cohomology sheaf
$ H^{j}(\pi\times id)_{*}\cB_{f'} $ satisfies the conditions of
the filtration
$ V $ in {\rm (2.1).}
Furthermore, if
$ \pi $ induces an isomorphism over
$ X \setminus D $, we have a canonical surjective morphism
$$
(H^{0}(\pi\times id)_{*}\cB_{f'};F,V) \to (\cB_{f};F,V),
\leqno(3.2.1)
$$
which induces isomorphisms for
$ p \in \bZ $,
$ \alpha > 0 $
}
$$
F_{p}V^{\alpha}H^{0}(\pi\times id)_{*}\cB_{f'} =
F_{p}V^{\alpha}\cB_{f}
\leqno(3.2.2)
$$

\medskip\noindent
{\it Proof.}
The first assertion follows from [23], 3.3.17.
Forgetting
$ V $, (3.2.1) is induced by the canonical morphism
$$
H^{0}\pi_{*}(\cO_{X'},F) \to (\cO_{X},F),
\leqno(3.2.3)
$$
whose restriction to
$ F_{-n} $ is given by the trace morphism
$ \pi_{\ssbull}\omega_{X'/X} \to \cO_{X} $.
The compatibility with
$ V $ follows from the functoriality of the filtration
$ V $, see [23], 3.1.5.
The kernel of (3.2.1) is supported on
$ D $ so that
$ V^{\alpha} $ of the kernel vanishes for
$ \alpha > 0 $, see [23], 3.1.3.
So we get (3.2.2) because morphisms of Hodge modules are
bistrictly compatible with
$ (F,V) $, see [23], 3.3.3--5.

\medskip\noindent
{\bf 3.3.~Remark.}
The bistrictness implies that
$ F_{p}V^{\alpha}H^{0} = H^{0}F_{p}V^{\alpha} $,
see [23], 1.2.13.
So the left-hand side of (3.2.2) for
$ p = -n = -m $ is
$ \pi_{\ssbull}(\omega_{X'/X}\otimes V^{\alpha}\cO_{X'}) $
by (3.1.4).
The bistrictness implies also
$$
R^{j}\pi_{\ssbull}(\omega_{X'/X}\otimes V^{\alpha}\cO_{X'})
= 0 \quad \text{for}\,\, j > 0.
\leqno(3.3.1)
$$
Here
$ V^{\alpha}\cO_{X'} $ can be replaced with
$ \cO_{X'}(-\msum_{i} [(\alpha -\varepsilon)m_{i}]D'_{i}) $
by (2.2--3).
This gives another proof of (1.1.3).

\medskip\noindent
{\bf 3.4~Proof of Theorem (0.1).}
By (2.3) and (3.2), we have for
$ \alpha > 0 $
$$
F_{-n}V^{\alpha}\cB_{f} = F_{-n}V^{\alpha}H^{0}
(\pi\times id)_{*}\cB_{f'} =
\pi_{\ssbull}(\omega_{X'/X}\otimes V^{\alpha}\cO_{X}).
\leqno(3.4.1)
$$
The left-hand side of (3.4.1) is
$ V^{\alpha}\cO_{X} $ by definition.
Using (1.1.1), (2.3.1), this implies
$$
V^{\alpha}\cO_{X} = \cJ((\alpha -\varepsilon)D)\quad
\text{for}\,\, \alpha \in \bQ\quad
\text{and}\,\, 0 < \varepsilon \ll 1.
$$
This completes the proof of Theorem (0.1).

\medskip\noindent
{\bf 3.5.~Proposition.}
{\it Assume
$ D $ reduced, and let
$ \tom_{D} = \rho_{*}\omega_{\tD} \,(\subset
\omega_{D}) $ as in the introduction.
Then
$$
\gathered
\tom_{D}
= \pi_{\ssbull}(W_{1}\omega_{X'}(D'_{\red}))/\omega_{X}
\subset \omega_{D}
= \omega_{X}(D)/\omega_{X},
\\
\omega_{D}/\tom_{D}
= \omega_{X}(D)/\pi_{\ssbull}(W_{1}\omega_{X'}(D'_{\red})).
\endgathered
\leqno(3.5.1)
$$
Here
$ \omega_{X}(D) = \omega_{X}\otimes \cO_{X}(D) $ and
$ W $ denotes the weight filtration on the logarithmic forms
$ \omega_{X'}(D'_{\red}) $ {\rm (}see {\rm [6])} so that
$ W_{1}\omega_{X'}(D'_{\red}) = (\partial f')\omega_{X'}
(D') $ with
$ (\partial f') $ the ideal generated by the partial
derivatives of
$ f' := \pi^{*}f $.
If the proper transform
$ D'' $ of
$ D $ is smooth {\rm (}i.e. if the proper transforms of the
irreducible components of
$ D $ do not intersect each other{\rm )}, then
}
$$
\gathered
\tom_{D}
= \pi_{\ssbull}(\omega_{X'}(D''))/\omega_{X}
\subset \omega_{D}
= \omega_{X}(D)/\omega_{X},
\\
\omega_{D}/\tom_{D}
= \omega_{X}(D)/\pi_{\ssbull}(\omega_{X'}(D'')).
\endgathered
\leqno(3.5.2)
$$

\medskip\noindent
{\it Proof.}
If
$ D'' $ is smooth, we may take
$ \tD = D'' $, and the vanishing of higher direct image
$ R^{1}\pi_{\ssbull}\omega_{X'} $ implies the isomorphisms
of (3.5.2).
(A similar formula was shown in [3].)
In this case we can also show
$$
\pi_{\ssbull}(W_{1}\omega_{X'}(D'_{\red})) = \pi_{\ssbull}
(\omega_{X'}(D'')),
\leqno(3.5.3)
$$
because
$ \pi_{\ssbull}(\omega_{D'_{j}}) = 0 $ for any irreducible
component
$ D'_{j} $ of
$ D'_{\red} $ such that
$ \dim \pi (D'_{j}) < \dim D'_{j} $.
(This is easy if
$ \pi $ is obtained by iterating blow-ups with smooth centers.)
Finally, we can verify that
$ \pi_{\ssbull}(W_{1}\omega_{X'}(D'_{\red})) $ is
independent of the embedded resolution
$ (X',D') $.
So the assertion follows.

\medskip\noindent
{\bf 3.6.~Corollary.}
{\it Choosing a defining equation
$ f $ of
$ D $, we have
$$
\Gr_{V}^{1}(\omega_{D}/\tom_{D})
= V^{1}\omega_{X}/\pi_{\ssbull}((\partial f)\omega_{X'}),
\leqno(3.6.1)
$$
where
$ V^{1}\omega_{X} = V^{1}\cO_{X}\otimes\omega_{X} $.
If
$ D'' $ is smooth, then
}
$$
\Gr_{V}^{1}(\omega_{D}/\tom_{D})
= V^{1}\omega_{X}/\pi_{\ssbull}(\omega_{X'}(D'' - D')).
\leqno(3.6.2)
$$

\medskip\noindent
{\bf 3.7.~Remark.}
The formula (3.6.2) is essentially due to [3].
The sheaf
$ \pi_{\ssbull}(\omega_{X'/X}(D'' - D')) $ is called the
adjoint ideal, and has been studied in [28] and [11].

\bigskip\bigskip
\centerline{{\bf 4. Localizations of Nearby Cycle Sheaves}}

\bigskip\noindent
{\bf 4.1.}
With the notation of (2.1) we define the nearby and
vanishing cycles of
$ (\cO_{X},F) $ and
$ (\cB_{f},F) $ by
$$
\aligned
\psi_{f,\be(-\alpha)}(\cO_{X},F) =
\psi_{t,\be(-\alpha)}(\cB_{f},F)
&:=
\Gr_{V}^{\alpha}(\cB_{f},F[1])\,\, (\alpha \in (0,1]),
\\
\varphi_{f,1}(\cO_{X},F) =
\varphi_{t,1}(\cB_{f},F)
&:=
\Gr_{V}^{0}(\cB_{f},F),
\endaligned
\leqno(4.1.1)
$$
and
$ \psi_{f}(\cO_{X},F) = \mopls_{0<\alpha \le 1}
\psi_{f,\be(-\alpha)}(\cO_{X},F) $ (similarly for
$ \cB_{f} $).
These are the underlying filtered
$ \cD_{X} $-modules of the nearby and vanishing cycles of mixed
Hodge modules
$ \psi_{f}(\bQ_{X}^{H}[n]) $,
$ \varphi_{f,1}(\bQ_{X}^{H}[n]) $, see [23], [24].
Indeed, these functors correspond to the nearby and vanishing
cycles functors [7] for constructible sheaves, which are shifted
by
$ -1 $, see [15], [20].
The last shift is needed for the stability of perverse sheaves
[1] by these functors.

Assume now that
$ D $ is a divisor with normal crossings.
Let
$ E, E' $ be reduced divisors on
$ X $ such that
$$
E \cup E' = D_{\red}\quad\text{and}\quad
\codim E \cap E' \ge 2.
\leqno(4.1.2)
$$
Set
$ E'' = E \cap E' $.
Let
$ i : E \to X, j : X \setminus E \to X $ denote the
inclusion morphisms, and similarly for
$ i', j' $ with
$ E $ replaced by
$ E' $.
For
$ \alpha \in (0,1] $,
let
$$
E(\alpha) = D(\alpha) \cap E
\leqno(4.1.3)
$$
with
$ D(\alpha) $ as in (2.2.5), and similarly for
$ E'(\alpha), E''(\alpha) $.
Then

\medskip\noindent
{\bf 4.2.~Theorem.}
{\it
With the above notation and the assumption,
there are canonical isomorphisms
}
$$
\aligned
i_{*}i^{*}\psi_{f}(\bQ_{X}^{H}[n])
&\simto
j'_{*}j^{\prime *}\psi_{f}(\bQ_{X}^{H}[n]),
\\
j'_{!}j^{\prime *}\psi_{f}(\bQ_{X}^{H}[n])
&\simto
i_{*}i^{!}\psi_{f}(\bQ_{X}^{H}[n]).
\endaligned
\leqno(4.2.1)
$$

\medskip\noindent
{\it Proof.}
Since we have canonical morphisms by the property of open
direct images, it is enough to show (4.2.1) for the underlying
perverse sheaves or
$ \cD $-modules, and the assertion is local.

Let
$ M = \psi_{f}\cO_{X} $.
Take local coordinates as in (2.2).
For
$ \mu = (\mu_{1}, \dots, \mu_{n}) \in \bQ^{n} $, let
$$
M^{\mu} = \mcap_{i}(\mcup_{j}
\Ker((\partial_{i}x_{i}-\mu_{i})^{j} : M \to M)).
$$
Then
$ M $ is generated by
$ \mopls_{\mu}M^{\mu} $, and is a completion of
$ \mopls_{\mu}M^{\mu} $ (i.e. regular singular of normal
crossing and quasi unipotent type, see [24], 3.2).
So we get an infinite direct sum decomposition by the
$ M^{\mu} $, which is compatible with the Hodge filtration
$ F $.

To simplify the argument, we first show the case
$ D $ reduced (i.e. semistable) and
$ f = x_{1} \cdots x_{n} $.
In this case
$ M^{\mu} = 0 $ unless
$ \mu \in \bZ^{n} $ (i.e.
$ M $ has unipotent local monodromies).
For a subset
$ I $ of
$ \{1, \dots, n\} $,
let
$ M_{I} = M^{\mu} $ with
$ \mu_{i} = 0 $ if
$ i \in I $ and
$ \mu_{i} = 1 $ otherwise.
For
$ i \in I $,
we have morphisms
$$
\partial_{i} : M_{I\setminus \{i\}} \to M_{I},\quad
x_{i} : M_{I} \to M_{I\setminus \{i\}}.
$$
It is well known that
$ \psi_{f}\cO_{X} $ is determined by
$ \{M_{I}\} $ with the above morphisms.
Using the calculation in (2.2), it is easy to see that
$$
M_{I} = \bC[s]/\bC[s]s^{|I|}.
\leqno(4.2.2)
$$
Furthermore,
$ \partial_{i} $ is given by the multiplication by
$ s $,
and
$ x_{i} $ by the projection.

Let
$ J, J' $ be the subsets of
$ \{1, \dots, n\} $ such that
$ E = \mcup_{i\in J} x_{i}^{-1} $,
$ E' = \mcup_{i\in J'} x_{i}^{-1} $ locally.
Let
$$
M(E_{!}) = j_{!}j^{*}\psi_{f}\cO_{X},\quad
M(E'_{*}) = j'_{*}j^{\prime *}\psi_{f}\cO_{X},
$$
and define
$ M(E_{!})_{I} $,
$ M(E'_{*})_{I} $ as above.
Here
$ j_{*} $ is the (regular singular) meromorphic direct image,
and
$ j_{!} = \bD j_{*}\bD $.
(Note that
$ i_{*}i^{!} $ and
$ i_{*}i^{*} $ in (4.2.1) are defined by the mapping cone of
$ j_{!}j^{*} \to id $ and the shifted mapping cone of
$ id \to j'_{*}j^{\prime *} $ respectively.)
Then
$$
M(E_{!})_{I} = \bC[s]s^{|I\cap J|}/\bC[s]s^{|I|},\quad
M(E'_{*})_{I} = \bC[s]/\bC[s]s^{|I\setminus J'|},
$$
so that we have the bijective morphisms
$$
\alignat 2
&\partial_{i} : M(E_{!})_{I} \simto
M(E_{!})_{I\cup \{i\}} \quad
&\text{for}\,\,
&i \in J \setminus I,
\\
&x_{i} : M(E'_{*})_{I\cup \{i\}} \simto
M(E'_{*})_{I} \quad
&\text{for}\,\,
&i \in J' \setminus I.
\endalignat
$$
Thus we get the exact sequence
$$
0 \to j_{!}j^{*}\psi_{f}\cO_{X} \to \psi_{f}\cO_{X} \to
j'_{*}j^{\prime *}\psi_{f}\cO_{X} \to 0,
\leqno(4.2.3)
$$
and the assertion follows in this case.

In general, we use
$ M_{I}^{\mu} $ instead of
$ M_{I} $ where
$ M_{I}^{\mu} = M^{\mu'} $ with
$ \mu' = \mu - 1_{I} $.
Here
$ \mu \in (0,1]^{n} $ and
$ I \subset \{1, \dots, n\} $ satisfy
$ \mu_{i} = 1 $ for
$ i \in I $, and
$ 1_{I} \in \bZ^{n} $ is defined by
$ (1_{I})_{i} = 1 $ if
$ i \in I $ and
$ 0 $ otherwise.
For
$ M = \psi_{f}\cO_{X} $, we have
$$
M_{I}^{\mu} = \bC[s]/\bC[s]s^{|I|} \quad
\text{if}\,\, \mu_{i} + m_{i}\alpha \in \bZ \,\,
\text{for any}\,\, i,
$$
and
$ M_{I}^{\mu} = 0 $ otherwise, see [24], 3.3.
Here we assume
$ m_{i} \ne 0 $ for any
$ i $.
Then the argument is similar.
This finishes the proof of (4.2).

\medskip\noindent
{\bf 4.3.~Theorem.}
{\it With the notation and the assumption of {\rm (4.1),}
there are canonical isomorphisms
}
$$
\aligned
F_{1-n}(i_{*}i^{*}\Gr_{V}^{\alpha}\cB_{f})
&= F_{1-n}\Gr_{V}^{\alpha}\cB_{f}\otimes_{\cO_{X}}
\cO_{E(\alpha)},
\\
F_{1-n}(i_{*}i^{!}\Gr_{V}^{\alpha}\cB_{f})
&= F_{1-n}\Gr_{V}^{\alpha}\cB_{f}\otimes_{\cO_{X}}
\cO_{E(\alpha)}(-E''(\alpha)).
\endaligned
\leqno(4.3.1)
$$

\medskip\noindent
{\it Proof.}
We first consider the case
$ D $ reduced as above.
Since the Hodge filtration in the normal crossing case is
compatible with the decomposition by the
$ M^{\mu} $, and the Hodge filtration on
$ M_{I} $ is given by the order of
$ s $,
we can verify the first isomorphism of (4.3.1) using (4.2.2).
This implies the second isomorphism by (4.2) because the
difference between
$ j'_{!}j^{\prime *}\psi_{f}\cO_{X} $ and
$ j'_{*}j^{\prime *}\psi_{f}\cO_{X} $ is given by the
multiplication by
$ \prod_{i\in J'} x_{i} $ (see also [24], 3.10 for the
difference between the two open direct images in the normal
crossing case).
The argument in the general case is similar.
This finishes the proof of (4.3).

\bigskip\bigskip\centerline{{\bf 5. Spectrum}}

\bigskip\noindent
{\bf 5.1.}
Let
$ X $ be a complex manifold of dimension
$ n $,
and
$ f $ be a holomorphic function on
$ X $.
For
$ x \in X $,
let
$ H^{j}(F_{x},\bQ) $ denote the vanishing cohomology at
$ x $ (i.e.
$ F_{x} $ denotes the ``fiber'' of the Milnor fibration
around
$ x) $.
It has a canonical mixed Hodge structure (see [22], [24]).
Using mixed Hodge modules, it is given by
$$
H^{j}(F_{x},\bQ) = H^{j-n+1}i_{x}^{*}\psi_{f}
(\bQ_{X}^{H}[n]),
\leqno(5.1.1)
$$
where
$ i_{x} : \{x\} \to X $ denotes the inclusion, and we put
$ F^{p} = F_{-p} $.
(See (4.1.1) for the underlying filtered
$ \cD $-module of
$ \psi_{f}(\bQ_{X}^{H}[n]) $.
Note that the functor
$ \psi $ for mixed Hodge modules is shifted by
$ -1 $.)
It has also the action of the semisimple part
$ T_{s} $ of the monodromy
$ T $.
Let
$ H^{j}(F_{x},\bC)_{\lambda} $ denote the
$ \lambda $-eigenspace.
Put
$ \be(\alpha) = \exp(2\pi i\alpha) $.

The spectrum
$ \Sp(f,x) $ is a fractional Laurent polynomial
$ \msum_{\alpha \in \bQ} n_{\alpha}t^{\alpha} $ such that
$$
n_{\alpha} = \msum_{j} (-1)^{j-n+1} \dim \Gr_{F}^{p}
\tH^{j}(F_{x},\bC)_{\be(-\alpha)}\quad
\text{with}\,\, p = [n - \alpha],
\leqno(5.1.2)
$$
where
$ \tH^{j}(F_{x},\bC) $ denotes the reduced
cohomology.
This normalization of the spectrum is different from the one in
[27] by the multiplication by
$ t $,
and coincides with the one in [4], [25].

\medskip\noindent
{\bf 5.2.~Proposition.}
{\it
$ n_{\alpha} = 0 $ if
$ \alpha \ge n $ or
$ \alpha \le 0 $.
}

\medskip\noindent
{\it Proof.}
The assertion for
$ \alpha \le 0 $ is equivalent to
$ \Gr_{F}^{p}\tH^{j}(F_{x},\bC) = 0 $ for
$ p \ge n $.
The functor
$ i_{x}^{*} $ is defined by using co-Cech complex
of open direct images, and preserves the property that
$ F_{p} = 0 $ for
$ p \le - n $.
So the assertion is clear for
$ \alpha \le 0 $.
(Note that the filtration is shifted by
$ 1 $ in (4.1.1).)

For
$ \alpha \ge n $ we have the duality isomorphisms
(see [24], 2.6)
$$
\aligned
\bD(\psi_{f}\bQ_{X}^{H}[n])
&= (\psi_{f}\bQ_{X}^{H}[n])(n-1),
\\
\bD(\varphi_{f,1}\bQ_{X}^{H}[n])
&= (\varphi_{f,1}\bQ_{X}^{H}[n])(n),
\endaligned
\leqno(5.2.1)
$$
because
$ \bD(\bQ_{X}^{H}[n]) = (\bQ_{X}^{H}[n])(n) $.
Here
$ \bD $ is the functor which associates the dual.
We first show
$ \Gr_{F}^{p}\tH^{j}(F_{x},\bC)_{1} = 0 $ for
$ p \le 0 $.
The second isomorphism of (5.2.1) implies
$$
\bD\tH^{j}(F_{x},\bQ)_{1} =
H^{-j+n-1}i_{x}^{!}\varphi_{f,1}(\bQ_{X}^{H}[n])(n)
\leqno(5.2.2)
$$
because
$ \tH^{j}(F_{x},\bQ)_{1} = H^{j-n+1}i_{x}^{*}
\varphi_{f,1}(\bQ_{X}^{H}[n]) $, and
$ \bD $ exchanges
$ i_{x}^{*} $ and
$ i_{x}^{!} $.
By (2.1.4) we have the vanishing
of
$ F_{p} $ on the underlying
$ \cD $-module of
$ \varphi_{f,1}(\bQ_{X}^{H}[n]) $ for
$ p < 1 - n $, and the functor
$ i_{x}^{!} $ preserves the property that
$ F_{p} = 0 $ for
$ p < 1 - n $ (i.e.
$ F^{p} = 0 $ for
$ p > n - 1 $).
So the assertion for
$ \alpha \in \bZ $ follows from (5.2.2).

Similarly, we can show that
$ \Gr_{F}^{p}H^{j}(F_{x},\bC) = 0 $ for
$ p < 0 $,
using the first isomorphism of (5.2.1).
This completes the proof of (5.2).

\medskip\noindent
{\bf 5.3.~Remarks.}
(i) In the isolated singularity case, (5.2) is due to Steenbrink
[26].

(ii) We can show in general
$$
\Gr_{F}^{p}H^{j}(F_{x},\bC) = 0
\quad\text{for}\,\, p > j.
\leqno(5.3.1)
$$
Indeed, the pull-back to a divisor with a defining equation
$ g $ (e.g. a coordinate function) is represented by
$ C(\can : \psi_{g,1} \to \varphi_{g,1}) $, and the Hodge
filtration is shifted by
$ 1 $ for
$ \psi $.

\medskip\noindent
{\bf 5.4.}
Let
$ \pi : (X',D') \to (X,D) $ be an embedded resolution
such that
$ E := \pi^{-1}(x)_{\red} $ is a divisor with normal
crossings.
Let
$ E' $ be the reduced divisor on
$ X' $ such that
$ D'_{\red} = E \cup E' $ and
$ E'' := E \cap E' $ has codimension
$ \ge 2 $.
Let
$ E(\alpha) = E \cap D'(\alpha) $,
and similarly for
$ E'(\alpha), E''(\alpha) $.
With the notation of (1.1) we define
$$
\aligned
K_{\alpha}
&= \omega_{X'/X}\otimes\cO_{X'}
(-\msum_{i} [(\alpha -\varepsilon)m_{i}]D'_{i}))
\otimes\cO_{E(\alpha)},
\\
K'_{\alpha}
&= \omega_{X'/X}\otimes\cO_{X'}
(-\msum_{i} [(\alpha -\varepsilon)m_{i}]D'_{i}))
\otimes\cO_{E(\alpha)}(-E''(\alpha)),
\endaligned
\leqno(5.4.1)
$$
for
$ 0 < \varepsilon \ll 1 $.
In the algebraic case, we have by [4]
$$
\cG(D,x,\alpha) = \Gamma (E(\alpha), K_{\alpha}).
\leqno(5.4.2)
$$

\medskip\noindent
{\bf 5.5.~Proof of Theorem (0.5).}
By (2.3) and (4.1.1) we have canonical isomorphisms
$$
F_{1-n}(\psi_{f,\be(-\alpha)}\cO_{X'}) =
F_{1-n}\Gr_{V}^{\alpha}\cB_{f'} =
\Gr_{V}^{\alpha}\cO_{X'}\,\,
(\alpha\in (0,1]).
\leqno(5.5.1)
$$
Let
$ E(\alpha), E'(\alpha), E''(\alpha) $ be as in (4.1),
and
$ i : E \to X' $ be the inclusion morphism.
Then by (2.2--3) and (4.3) we get
$$
\aligned
K_{\alpha}
&= \omega_{X'/X}\otimes F_{1-n}\Gr_{V}^{\alpha}
\cB_{f'}\otimes \cO_{E(\alpha)}
\\
&\qquad
= \omega_{X'/X}\otimes F_{1-n}(i_{*}i^{*}\Gr_{V}^{\alpha}
\cB_{f'}),
\\
K'_{\alpha}
&= \omega_{X'/X}\otimes F_{1-n}\Gr_{V}^{\alpha}
\cB_{f'}\otimes \cO_{E(\alpha)}(-E''(\alpha))
\\
&\qquad
= \omega_{X'/X}\otimes F_{1-n}(i_{*}i^{!}\Gr_{V}^{\alpha}
\cB_{f'}).
\endaligned
$$
Taking the cohomology of these coherent sheaves, we get
$ F_{1-n} $ of the direct image of
$ i_{*}i^{*}\Gr_{V}^{\alpha}\cB_{f'} $,
$ i_{*}i^{!}\Gr_{V}^{\alpha}\cB_{f'} $ under
$ \pi : X' \to X $ by definition of the direct image of filtered
$ \cD $-modules (see [23] and also (3.2) above), because
$ F_{p} = 0 $ for
$ p < 1 - n $.
We have
$$
(i_{x})_{*}i_{x}^{*}\psi_{f}(\bQ_{X}^{H}[n]) =
(i_{x})_{*}i_{x}^{*}\pi_{*}\psi_{f'}(\bQ_{X'}^{H}[n]) =
\pi_{*}i_{*}i^{*}\psi_{f'}(\bQ_{X'}^{H}[n]),
$$
and similarly for
$ i_{x}^{!} $,
$ i^{!} $.
So we get the assertion.

\bigskip\bigskip
\centerline{{\bf References}}

\bigskip

\item{[1]}
Beilinson, A., Bernstein, J. and Deligne, P., Faisceaux pervers,
Ast\'erisque 100, Soc. Math. France, Paris, 1982.

\item{[2]}
Borel, A. et al., Algebraic
$ {\Cal D} $-modules, Perspectives in Math. 2, Academic Press,
1987.

\item{[3]}
Budur, N., thesis.

\item{[4]}
\SameAuthor, On Hodge spectrum and multiplier ideals, to appear in
Math. Ann.

\item{[5]}
\SameAuthor, Multiplier ideals and filtered
$ D $-modules, preprint.

\item{[6]}
Deligne, P., Th\'eorie de Hodge I, Actes Congr\`es Intern.
Math., Part 1 (1970), 425--430; II, Publ. Math. IHES, 40 (1971),
5--58; III, ibid. 44 (1974), 5--77.

\item{[7]}
\SameAuthor, Le formalisme des cycles \'evanescents, in SGA7 XIII
and XIV, Lect. Notes in Math. 340, Springer, Berlin, 1973, pp.
82--115 and 116--164.

\item{[8]}
Demailly, J.-P.,
$ L^{2} $-Methods and effective results in algebraic geometry,
Proc. Intern. Congr. Math. (Z\"urich, 1994), Birkh\"auser, Basel,
1995, pp. 817--827.

\item{[9]}
Denef, J. and Loeser, F., Motivic Igusa zeta functions, J. Alg.
Geom. 7 (1998), 505--537.

\item{[10]}
Ein, L., Multiplier ideals, vanishing theorems and applications,
Proc. Symp. Pure Math., A.M.S. 62 Part 1, (1997), 203--219.

\item{[11]}
Ein, L. and Lazarsfeld, R., Singularities of theta divisors
and the birational geometry of irregular varieties,
J. Amer. Math. Soc. 10 (1997), 243--258.

\item{[12]}
Ein, L., Lazarsfeld, R., Smith, K.E. and Varolin, D., Jumping
coefficients of multiplier ideals, preprint.

\item{[13]}
Grauert, H. and Riemenschneider, O., Verschwindungss\"atze
f\"ur analytische Kohomologiegruppen auf Komplexen R\"aumen,
Inv. Math. 11 (1970), 263--292.

\item{[14]}
Kashiwara, M.,
$ B $-functions and holonomic systems, Inv. Math. 38 (1976/77),
33--53.

\item{[15]}
\SameAuthor, Vanishing cycle sheaves and holonomic systems of
differential equations, Algebraic geometry (Tokyo/Kyoto, 1982),
Lect. Notes in Math. 1016, Springer, Berlin, 1983, pp. 134--142.

\item{[16]}
Kashiwara, M. and Kawai, T., Second-microlocalization and
asymptotic expansions, Proc. Internat. Colloq., Centre Phys.
(Les Houches, 1979), Lecture Notes in Phys., 126, Springer,
Berlin, 1980, pp. 21--76.

\item{[17]}
Koll\'ar, J., Singularities of pairs, Proc. Symp. Pure Math.,
A.M.S. 62 Part 1, (1997), 221--287.

\item{[18]}
Loeser, F., Quelques cons\'equences locales de la th\'eorie
de Hodge, Ann. Inst. Fourier 35 (1985) 75--92.

\item{[19]}
Malgrange, B., Le polyn\^ome de Bernstein d'une
singularit\'e isol\'ee, in Lect. Notes in Math. 459, Springer,
Berlin, 1975, pp. 98--119.

\item{[20]}
\SameAuthor, Polyn\^ome de Bernstein-Sato et cohomologie
\'evanescente, Analysis and topology on singular spaces, II,
III (Luminy, 1981), Ast\'erisque 101--102 (1983), 243--267.

\item{[21]}
Nadel, A.M., Multiplier ideal sheaves and K\"ahler-Einstein
metrics of positive scalar curvature, Ann. Math. 132 (1990),
549--596.

\item{[22]}
Navarro Aznar, V., Sur la th\'eorie de Hodge-Deligne, Inv.
Math. 90 (1987), 11--76.

\item{[23]}
Saito, M., Modules de Hodge polarisables, Publ. RIMS, Kyoto
Univ. 24 (1988), 849--995.

\item{[24]}
\SameAuthor, Mixed Hodge modules, Publ. RIMS, Kyoto Univ. 26
(1990), 221--333.

\item{[25]}
\SameAuthor, On
$ b $-function, spectrum and rational singularity,
Math. Ann. 295 (1993), 51--74.

\item{[26]}
Steenbrink, J.H.M., Mixed Hodge structure on the vanishing
cohomology, in Real and Complex Singularities (Proc. Nordic
Summer School, Oslo, 1976) Alphen a/d Rijn: Sijthoff \& Noordhoff
1977, pp. 525--563.

\item{[27]}
\SameAuthor, The spectrum of hypersurface singularity,
Ast\'erisque 179--180 (1989), 163--184.

\item{[28]}
Vaqui\'e, M., Irr\'egularit\'e des rev\^etements cycliques,
in Singularities (Lille, 1991), London Math. Soc. Lecture
Note Ser., 201, Cambridge Univ. Press, Cambridge, 1994,
pp. 383--419.

\item{[29]}
Varchenko, A., Asymptotic Hodge structure in the vanishing
cohomology, Math. USSR Izv. 18 (1982), 465--512.

\bigskip

Nero Budur

Department of Mathematics, University of Illinois at Chicago

851 South Morgan Street (M/C 249) Chicago, IL 60607-7045, USA

E-Mail: nero\@math.uic.edu

\bigskip

Morihiko Saito

RIMS Kyoto University, Kyoto 606-8502 Japan

E-Mail: msaito\@kurims.kyoto-u.ac.jp

\bigskip

\ver

\bye